\documentclass[11pt,letterpaper]{amsart}
\usepackage{preamble}

\usetikzlibrary{calc}
\definecolor{cycleblue}{HTML}{2463A0}
\definecolor{centerpurple}{HTML}{7852A0}
\definecolor{helperorange}{HTML}{D98519}
\definecolor{friendred}{HTML}{BB3946}
\definecolor{ink}{HTML}{25313C}

\counterwithin{equation}{section}

\begin{document}

\title[Tight Probabilistic Aspects of Friends-and-Strangers Graphs]{Tight Probabilistic Aspects of Friends-and-Strangers Graphs}
\author[Neil Krishnan]{Neil Krishnan}
\date{\today}

\maketitle

\begin{abstract}
Friends-and-strangers graphs, introduced by Defant and Kravitz, model arrangements of people on positions. We consider two graphs $X$ and $Y$ on $n$ vertices. The vertices of $X$ represent positions with edges denoting adjacencies while the vertices of $Y$ represent people with edges denoting friendships. Their \emph{friends-and-strangers graph} is denoted by $\mathsf{FS}(X,Y)$. Its vertices are all one-to-one placements of people on positions, and two placements are adjacent if one can swap two friends on adjacent positions to get from one placement to the other. 

A fundamental question is whether or not a friends-and-strangers graph is connected. A natural probabilistic analog is ``If $X$ and $Y$ are independently sampled Erd\H{o}s--R\'{e}nyi graphs with $n$ vertices and edge probability $p,$ what is the threshold probability for $\mathsf{FS}(X,Y)$ to be connected with high probability?'' Alon, Defant, and Kravitz (2023) showed that $p = n^{-1/2+o(1)}$ is the threshold probability and conjectured that $p = \Theta(n^{-1/2}).$ Milojevi\'{c} (2024) disproved this by showing $p = \Omega(n^{-1/2}\sqrt{\log n})$ is necessary for connectivity with high probability. In this work, we prove a matching upper bound, thus establishing that the threshold is $p = \Theta(n^{-1/2}\sqrt{\log n}).$ 
\end{abstract}

\section{Introduction}
Friends-and-strangers graphs were formalized by Defant and Kravitz \cite{defant2021friends}.
\begin{definition}
    Given two graphs $X$ and $Y$ each on $n$ vertices, the \emph{friends-and-strangers graph} $\FS(X,Y)$ is defined as follows.
    Its vertex set is the set of bijections $\sigma:V(X)\to V(Y)$, and $\sigma$ and $\tau$ are adjacent if and only if $\sigma=\tau\circ(i'\,j')$ for some $i',j'\in V(X)$ with $(i',j')\in E(X)$ and $(\sigma(i'),\sigma(j'))\in E(Y)$.
    In this case, we say $\sigma$ and $\tau$ differ by an \emph{$(X,Y)$-friendly swap}.
\end{definition}
Throughout this paper, we will adopt two conventions. First, we identify both vertex sets with $[n]$, i.e., $V(X) = V(Y) = [n].$ Thus, the bijections will actually be permutations on $[n].$ Second, when referring to vertices of $X$ we will add a prime, i.e., $v' \in V(X),$ and when referring to vertices of $Y$ we will omit the prime, i.e., $v \in V(Y).$

\begin{remark}
Friends-and-strangers graphs have a natural interpretation. Let the graph $X$ represent $n$ distinct positions, and let its edges denote adjacencies between positions. Similarly, let the graph $Y$ represent $n$ people, and let edges denote friendships. Friendship is not transitive in our problem. A bijection from $V(X)$ to $V(Y)$ is a placement of people on the positions. Two placements are an $(X,Y)$-friendly swap apart if two friends standing on adjacent positions can swap positions to get from one placement to the other. The graph $\FS(X,Y)$ thus has a vertex set consisting of the placements and the edges are all of the pairs of placements differing by an $(X,Y)$-friendly swap. We will use this interpretation throughout this paper.
\end{remark}

\begin{example}
Let $X = Y =P_3$ be the path graph on $3$ vertices. Then $\FS(X,Y)$ is the graph shown in \cref{fig:FSex}. This example was also mentioned in Defant and Kravitz \cite[Example 1.1]{defant2021friends}.
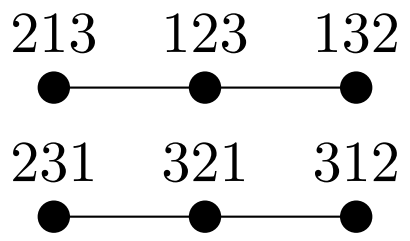
\begin{figure}[htbp]
    \centering
    \begin{tikzpicture}
        \draw (0,0) node[above,yshift=0.1 cm]{$213$}--(1,0) node[above,yshift=0.1 cm]{$123$}--(2,0) node[above,yshift=0.1 cm]{$132$};
        \filldraw (0,0) circle(0.1 cm) (1,0) circle(0.1 cm) (2,0) circle(0.1 cm);
    \end{tikzpicture}
    
    \vspace{0.1 cm}
    
    \begin{tikzpicture}
        \draw (0,0) node[above,yshift=0.1 cm]{$231$}--(1,0) node[above,yshift=0.1 cm]{$321$}--(2,0) node[above,yshift=0.1 cm]{$312$};
        \filldraw (0,0) circle(0.1 cm) (1,0) circle(0.1 cm) (2,0) circle(0.1 cm);
    \end{tikzpicture}

    \caption{The graph $\FS(P_3,P_3)$. The labels on the vertices mark the bijection from the positions to the people the vertex corresponds to, e.g., $312$ represents person $3$ is at position $1,$ person $1$ is at position $2,$ and person $2$ is at position $3.$}
    \label{fig:FSex}
\end{figure}
\end{example}

% Many aspects of friends-and-strangers graphs have been examined.
% \Neil{Diameter,Girth,Connectivity,Extremal \href{https://share.gemini.google/n5ghfMNeAHlR}{gemini}}

In their initial paper, Defant and Kravitz \cite{defant2021friends} studied the basic properties of $\FS(X,Y)$ and characterized the connected components of $\FS(P_n,Y)$ and $\FS(C_n,Y)$ using acyclic and toric acyclic orientations. Bangachev \cite{bangachev2022asymmetric} investigated extremal minimum-degree conditions guaranteeing connectivity, while Zhao, You, Lin, and Zhang \cite{zhao2025connectedness} subsequently established sufficient conditions for higher connectivity in terms of the maximum degree of $X$ and the vertex connectivity of $Y$. Krishnan and Li \cite{krishnan2026connectivity} studied the vertex connectivity of friends-and-strangers graphs, obtaining tight bounds for certain classes of friends-and-strangers graphs as well as sufficient extremal and probabilistic conditions. Regarding diameter, Jeong \cite{jeong2022diameters} constructed friends-and-strangers graphs with connected components of exponentially large diameter, while Akella and Li \cite{akella2026diameter} established polynomial diameter bounds under suitable minimum-degree assumptions. More recently, Mawhinney \cite{mawhinney2026connectivity} resolved a conjecture of Bangachev, proving that $\FS(X,Y)$ is connected whenever $X,Y$ are connected and
\[
2\min\{\delta(X),\delta(Y)\}+3\max\{\delta(X),\delta(Y)\}\ge 3n,
\]
and obtained asymptotically sharp diameter bounds under the same assumptions. Finally, Yang and Zhang \cite{yang2025pspace} showed that deciding whether two configurations belong to the same connected component of $\FS(X,Y)$ is PSPACE-complete.

\subsection{Previous Probabilistic Results}
Most relevant to our paper are previous probabilistic results on friends-and-strangers graphs. Let $\G(H,p)$ be a random subgraph of $H$ constructed by including each edge independently with probability $p$. We write $\G(n,p)=\G(K_n,p)$. The original probabilistic question is the following.
\begin{question}[Canonical Question]
    Let $X,Y$ be independently sampled from $\G(n,p)$ where $p = p(n).$ What is the threshold probability $p$ for $\FS(X,Y)$ to be connected?
\end{question}
In 2021, Alon, Defant, and Kravitz \cite{alon2021extremal} proved the following.
\begin{theorem}[\protect{\cite[Theorem 1.1]{alon2021extremal}}]\label{alonconnectivitywhp}
    Fix some small $\varepsilon > 0$. Let $X$ and $Y$ be independently-chosen random graphs in $\G(n,p),$ where $p=p(n)$ depends on $n.$ If
    $$p\leq \frac{2^{-1/2}-\varepsilon}{n^{1/2}},$$
    then $\FS(X,Y)$ is disconnected with high probability. If
    $$p\geq \frac{\exp(2(\log n)^{2/3})}{n^{1/2}},$$
    then $\FS(X,Y)$ is connected with high probability.
\end{theorem}
Notice that this implies the threshold probability is $n^{-1/2 + o(1)}.$ Alon, Defant, and Kravitz conjectured that the true threshold was $\Theta(n^{-1/2})$ but this was disproved in 2023 by Milojevi\'{c} \cite{milojevic2022connectivity}.
\begin{theorem} [\protect{\cite[Theorem 1.2]{milojevic2022connectivity}}] \label{milojevicdisconnectedwhp}
    There exists a constant $\varepsilon > 0$ with the following property. For a large positive integer $n$ and
    $$p<\varepsilon \left( \frac{\log n}{n} \right)^{1/2},$$
    if we choose random graphs $X$ and $Y$ in $\G(n,p)$ independently, $\FS(X,Y)$ has an isolated vertex with high probability.
\end{theorem}
Milojevi\'{c} further conjectured that $n^{-1/2} \log^{1/2} n$ was the correct probability threshold \cite[Conjecture 6.1]{milojevic2022connectivity}.

An important variation of the canonical question is the \emph{asymmetric version}, sampling $X \sim \G(n,p)$ and $Y \sim \G(n,q)$ independently.
% We can also consider the combination of the two where $X \sim G(K_{r,r},q)$ and $Y \sim G(K_{r,r},q)$ for $p = p(r)$ and $q = q(r)$ which is the \emph{asymmetric bipartite version} of the canonical question.
% In the same 2023 paper, Milojevi\'{c} considered the bipartite version of the canonical question.
% \begin{theorem}[\protect{\cite[Theorem 1.3]{milojevic2022connectivity}}]\label{milojevicbipartiteconnwhp}
%     Let $X$ and $Y$ be independently-chosen random graphs in $\mathcal{G}(K_{r,r},p),$ where $p=p(r)$ depends on $r.$ If
%     $$p\geq \frac{\exp(10(\log n)^{4/5})}{n^{1/2}},$$
%     then $\FS(X,Y)$ has two connected components with high probability.
% \end{theorem}
In 2022, Wang and Chen considered this asymmetric version.

\begin{theorem} [\protect{\cite[Theorem 1.3]{wang2023connectivity}}] \label{wangchenasymconnectedwhp}
Fix some small $\varepsilon > 0,$ and let $X$ and $Y$ be independently chosen random graphs in $\mathcal{G}(n,p_1)$ and $\mathcal{G}(n,p_2),$ respectively, where $p_1 = p_1(n)$ and $p_2=p_2(n)$ depend on $n$. Let $p_0 = \exp(2(\log n)^{2/3})n^{-1/2}.$ If either
$$p_1p_2 \leq \frac{(1-\varepsilon)/2}{n} \text{ and } p_1,p_2 \gg \frac{\log n}{n},$$
or
$$\min\{p_1,p_2\} \leq \frac{\log n + c(n)}{n} \text{ for some }c(n) \rightarrow -\infty,$$
then $\FS(X,Y)$ is disconnected with high probability. If $$p_1p_2 \geq p_0^2 \quad \text{and} \quad p_1,p_2\geq \frac{1}{\ell}p_0,$$ where $\ell = \frac{(\log n)^{1/3}}{2},$ then $\FS(X,Y)$ is connected with high probability.
\end{theorem}

\subsection{Heuristic Argument}
A heuristic argument shows why $pq=\Theta(\log n/n)$ is a reasonable threshold. We compare
$\FS(X,Y)$ with a random graph having the same numbers of vertices
and edges. $\FS(X,Y)$ has $n!$ vertices and exactly
$$|E(X)|\,|E(Y)|(n-2)!$$ edges. Indeed, to obtain an edge, we choose a position
edge $(u',v')\in E(X)$, a friendship $(a,b)\in E(Y)$, and a placement of
the remaining $n-2$ people on the remaining positions. The two ways
of placing $a,b$ on $u',v'$ give the endpoints of one friendly swap edge,
and every edge arises uniquely in this way. 

Thus, for independent
$X\sim\G(n,p)$ and $Y\sim\G(n,q)$, the expected average degree of
$\FS(X,Y)$ is $d=\binom{n}{2}pq$. A uniformly random graph on $n!$
vertices becomes connected around average degree $\log(n!)$
\cite{bollobas2001random}. Equating $\binom{n}{2}pq$ with
$\log(n!)=\Theta(n\log n)$ gives precisely $pq=\Theta(\log n/n)$.

\subsection{Main Result}
In this paper, we fully resolve the asymmetric upper bound.

\begin{theorem}\label{thm:main}
There exists a constant $C>0$ such that the following holds.
Let $X\sim\G(n,p)$ and $Y\sim\G(n,q)$ be independently
sampled Erd\H{o}s--R\'enyi random graphs. If
\[
pq\ge C\frac{\log n}{n},
\]
then $\FS(X,Y)$ is connected with high probability.
\end{theorem}

This is tight because $\FS(X,Y)$ is disconnected if it has an isolated vertex which occurs if there is a graph packing of $X$ and $Y$ into $K_n$ such that neither share an edge \cite{alon2021extremal}. The following mainly follows from Bollob\' as, Janson, and Scott \cite{bollobas2016packing} who determine the probability of such a packing existing when $p$ and $q$ are bounded away from $1$ (leaving the remaining case where one of $p$ or $q$ is close to $1$).

\begin{proposition}[\protect{\cite[Corollary 1.3]{baldissera2026connectivity}}]
There exists a constant $C'>0$ such that the following holds.
Let $X\sim\G(n,p)$ and $Y\sim\G(n,q)$ be independently
sampled Erd\H{o}s--R\'enyi random graphs. If
\[
pq\le C'\frac{\log n}{n},
\]
then $\FS(X,Y)$ is disconnected with high probability.
\end{proposition}

\begin{remark}\label{rem:diff}
    Recently, Baldissera, Dias do Couto, Griffiths, and Moreira \cite{baldissera2026connectivity} concurrently established the same connectivity threshold $pq=\Theta(\log n/n)$. Their proof extends the embedding-based approach of Alon, Defant, and Kravitz through a technically involved analysis of expansion in induced subgraphs of friendship neighborhoods, which is then used to embed auxiliary graphs.

By contrast, our proof is considerably more elementary than both their argument and the original approach of Alon, Defant, and Kravitz. We introduce a technique called \emph{movability}, which reduces connectivity to the ability to move individual people between adjacent positions. This allows us to replace the elaborate embedding constructions of \cite{alon2021extremal,baldissera2026connectivity} with two simple arms joined by a friendship edge. In particular, our argument avoids Janson's inequality and relies only on standard concentration bounds, elementary random graph properties, and union bounds.

In terms of similarities, both proofs use a partition to
obtain independence when testing friendship neighborhoods. They use
expansion in these neighborhoods to embed auxiliary graphs. With our
movability reduction we construct our auxiliary graphs directly using vertex connectivity.
\end{remark}

\section{Preliminaries}\label{sec:preliminaries}
All graphs are finite and simple.
We use the following standard graph-theoretic notation:
\begin{itemize}
    \item For a graph $G$ and $v\in V(G)$, let $N_G(v)$
    denote the set of neighbors of $v$, $d_G(v)$ denote the degree of $v$, and
    $\delta(G)$ denote the minimum degree of $G$.

    \item For $U\subseteq V(G)$, let $G[U]$ denote the
    subgraph of $G$ induced by $U$.

    \item A graph $G$ is \emph{$k$-vertex-connected} if
    $|V(G)|>k$ and deleting any set of fewer than
    $k$ vertices leaves a connected graph.
    Its \emph{vertex connectivity} $\kappa(G)$ is the
    largest $k$ for which $G$ is $k$-vertex-connected.
    \item An event holds \emph{with high
probability} if its probability tends to $1$ as $n\to\infty$.
\end{itemize}

All constants in this paper are independent of $n$; a sufficiently large
constant may be increased to satisfy finitely many requirements.

In the remainder of this section, we provide some probabilistic lemmas that will be used in later proofs. All of these statements are standard, and variations of these proofs may be found in \cite{bollobas2001random}; we include the quantitative versions needed below.

The main concentration inequality we use is the Chernoff bound.

\begin{lemma}[Chernoff bounds]\label{lem:chernoff}
Let $Z$ be a binomial or hypergeometric random variable with mean
$\mu$. For every $0<\varepsilon<1$,
\[
 \mathbb P\bigl(Z\le(1-\varepsilon)\mu\bigr)
 \le e^{-\varepsilon^2\mu/2},\qquad
 \mathbb P\bigl(Z\ge(1+\varepsilon)\mu\bigr)
 \le e^{-\varepsilon^2\mu/3}.
\]
\end{lemma}

The following lemmas give conditions for the following properties to hold:
\begin{itemize}
    \item Lemma~\ref{lem:random-vertex-connectivity}: $k$-vertex-connectivity,
    \item Lemma~\ref{lem:random-common-neighbors}: every pair of vertices having at least $k$ common neighbors,
    \item Lemma~\ref{lem:random-set-edges}: every pair of large disjoint sets of vertices connected by an edge,
    \item Lemma~\ref{lem:random-even-paths}: every pair of vertices connected by a short path of even length.
\end{itemize}

\begin{lemma}\label{lem:random-vertex-connectivity}
For every fixed $A>0$, there exists a constant $C_A>0$ with the
following property. Let $1\le k=k(n)\le\log n$, and let
$G\sim\G(n,p)$ where $p=p(n)\ge C_A n^{-1}\log n$.
Then $G$ is $k$-vertex-connected with probability at least $1-n^{-A}$.
\end{lemma}
\begin{proof}
For sufficiently large $C_A$, Chernoff and a union bound give
$\mathbb P(\delta(G)<2k)\le n e^{-\Omega(np)}$ as $k\leq \log n.$ 
Suppose $\delta(G)\ge2k$ but deleting a set $D$, $|D|<k$,
disconnects $G$. A smallest remaining component, with vertex set
$T$, satisfies that each vertex of $T$ has degree at least $k$ in $G[T]$. Since $T$ is a smallest remaining component,
$k\le t:=|T|\le(n-|D|)/2$.
For fixed $D,T$, the probability no edge joins $T$ to $V(G)\setminus(D\cup T)$ is at most $(1-p)^{t(n/3)} \leq e^{-npt/3}$ for large $n$ since $|V(G)\setminus(D\cup T)| \geq n/3$ for large $n$.
There are at most $kn^{k+t}$ choices of $D,T$ of these sizes.
Consequently, the failure probability is at most
\[
 n e^{-\Omega(np)}+
 k\sum_{t=k}^{\lfloor n/2\rfloor}n^{k+t}e^{-npt/3}
 \le n e^{-\Omega(np)}+
 k\sum_{t=k}^{\lfloor n/2\rfloor}n^{2t}e^{-npt/3}.
\]
Taking $C_A$ sufficiently large makes this at most $n^{-A}$.
\end{proof}

\begin{lemma}\label{lem:random-common-neighbors}
There exists a constant $C>0$ with the following property.
Let $k=k(n)\ge1$, and let $G\sim\G(n,p)$ where
\[
 p=p(n)\ge Cn^{-1/2}\sqrt{k+\log n}.
\]
Then every two distinct vertices of $G$ have at least $k$
common neighbors with high probability.
\end{lemma}
\begin{proof}
For fixed distinct vertices $u,v$, the number of their common
neighbors has distribution $\operatorname{Bin}(n-2,p^2)$, with mean
$\mu=(n-2)p^2$. For sufficiently large $C$ and $n$, we have
$\mu\ge 2k$ and $\mu\ge 24\log n$. By Chernoff,
\[
 \mathbb P\bigl(|N_G(u)\cap N_G(v)|<k\bigr)
 \le e^{-\mu/8}\le n^{-3}.
\]
A union bound over all pairs of vertices completes the proof.
\end{proof}

\begin{lemma}\label{lem:random-set-edges}
Let $1\le t\le n/2$ be an integer, and let $G\sim\G(n,p)$ where
$p\ge 3t^{-1}\log n$. Then, with probability at least $1-n^{-t}$,
every two disjoint sets of $t$ vertices have an edge between them.
\end{lemma}

\begin{proof}
For fixed disjoint sets, the probability of no edge between them is
$(1-p)^{t^2}\le e^{-pt^2}$. A union bound gives failure at most
$n^{2t}e^{-pt^2}\le n^{-t}$.
\end{proof}

\begin{lemma}\label{lem:random-even-paths}
There exists a constant $C>0$ with the following property.
Let $G\sim\G(n,p)$ where $p=p(n)\ge Cn^{-1}\log n$.
Then, with high probability, every two distinct vertices of $G$
are joined by a simple path of even length at most
\[
 2\left\lceil\frac{2\log n}{\log\log n}\right\rceil.
\]
\end{lemma}
\begin{proof}
Set $r=\lceil 2\log n/\log\log n\rceil$. Let $V(G)=A\sqcup B\sqcup J$ be an ordered balanced partition, i.e. with part sizes differing by
at most one. Let such a partition be \emph{good} if
\[
 \operatorname{diam}G[A],\operatorname{diam}G[B]\le r-1
\]
and there is a triangle $xyz$ with $x\in A$, $y\in B$, and $z\in J$.

For a good partition, every $u\in A$ and $v\in B$ are joined
by a simple even-length path of length at most $2r$. Indeed,
take shortest paths from $u$ to $x$ in $G[A]$ and from $y$ to
$v$ in $G[B]$. Joining these paths through $xy$ or through
$xzy$ gives two simple $u$--$v$ paths of opposite parity, each
of length at most $2(r-1)+2=2r$. Thus, it suffices to show every pair of vertices $u$ and $v$ satisfies $u \in A$ and $v \in B$ for some good partition $A\sqcup B\sqcup J.$

Because $p \geq Cn^{-1} \log n,$ for a fixed balanced partition, the diameter bound
\cite[Theorem 10.18]{bollobas2001random} gives $\operatorname{diam}G[A],\operatorname{diam}G[B]\le r-1$ with high probability. Although Theorem~10.18 is stated for the fixed-edge model (i.e., their graph is randomly chosen from graphs with exactly $M$ edges), we can use a Chernoff bound to condition on having at least $M$ edges in induced subgraphs of $G.$

A standard second-moment
calculation gives the required triangle with high probability,
since $np\ge C\log n$. Thus each fixed partition is good with
probability $1-o(1)$, uniformly over balanced partitions. By
linearity of expectation, the expected fraction of partitions
that fail to be good is $o(1)$. Markov's inequality then shows
that this fraction is less than $1/10$ with high probability.

Fix a $G$ where the fraction of bad balanced partitions is at most $1/10$. For any distinct vertices $u,v$,
a fraction $1/9+o(1)$ of the ordered balanced partitions have
$u\in A$ and $v\in B$. Since this exceeds $1/10$ for sufficiently
large $n$, at least one of those partitions is good.
\end{proof}

\section{Movability}
A critical idea in friends-and-strangers graphs is exchangeability which was first defined by Alon, Defant, and Kravitz \cite{alon2021extremal}. Let $X$ and $Y$ be two graphs on $n$ vertices and fix a bijection $\sigma:V(X) \to V(Y).$ Vertices $u,v \in V(Y)$ with $(\sigma^{-1}(u),\sigma^{-1}(v)) \in E(X)$ are \textit{$(X,Y)$-exchangeable from $\sigma$} if there is a sequence of $(X,Y)$-friendly swaps transforming $\sigma$ into $(u\,v)\circ\sigma$.
\begin{lemma} [\protect{\cite[Lemma 2.9]{alon2021extremal}}]\label{lem:exchangeable}
    Let $X$ and $Y$ be two graphs on $n$ vertices where $X$ is connected. Assume that for all vertices $u,v \in V(Y),$ and all permutations, $\sigma,$ where $(\sigma^{-1}(u),\sigma^{-1}(v)) \in E(X),$ we have that $u$ and $v$ are $(X,Y)$-exchangeable from $\sigma$. Then $\FS(X,Y)$ is connected.
\end{lemma}
This notion is powerful because it reduces the global problem of connectivity into a more localized problem of swapping just two people on adjacent positions in a placement.
Movability is an even more local condition building on exchangeability.

\begin{definition}
    Let $X$ and $Y$ be two graphs on $n$ vertices. Fix a bijection $\sigma: V(X) \to V(Y).$ Let $a,b \in V(Y)$ be distinct vertices, write $a'=\sigma^{-1}(a)$ and $b'=\sigma^{-1}(b)$, and suppose $(a',b')\in E(X)$. Then, we write \emph{$a \to b'$ is $(X,Y,k)$-movable from $\sigma$} if for every set $S \subseteq V(Y)\setminus\{a,b\}$ of size at most $k$ there exists a series of friendly swaps moving $a$ from $a'$ to $b'$ which never swaps with any person in $S.$ We will refer to $S$ as the people we \emph{protect}.
\end{definition}
In what follows $k$ will be relatively small, in particular, $k \leq \log n.$ Thus, for sufficiently well connected graphs, movability is a substantially weaker requirement than exchangeability as we do not have to maintain the positions of the other people. Under restrictions (which are light from a probabilistic standpoint), we can show that movability implies exchangeability and thus, applying Lemma~\ref{lem:exchangeable}, movability implies $\FS(X,Y)$ is connected.

The idea behind the following lemmas is to first have a particular pair $G$ and $H$ with $u,v \in V(H)$ of friends-and-strangers graphs and a bijection $\rho: V(G) \to V(H),$ so that $(u\,v)\circ \rho$ and $\rho$ are connected in $\FS(G,H).$ The idea then is to use movability to get into a situation where we can embed $G$ and $H$ into $X$ and $Y$ respectively, swap $u$ and $v$, and then reverse the movability steps to just get the original bijection with $u$ and $v$ swapped, implying exchangeability.

\begin{lemma}\label{lem:odd-cycle-exchange}
    Let $r\geq 1$, and let $G$ and $H$ be graphs on $2r+1$ vertices where $G$ is a cycle and $H$ has two vertices $u$ and $v$ whose neighborhoods are both the other $2r-1$ vertices. Let $\rho: V(G) \to V(H)$ be a bijection where $\rho^{-1}(u)$ and $\rho^{-1}(v)$ are adjacent. Then, there exist friendly swaps in $\FS(G,H)$ turning $\rho$ into $(u\,v) \circ \rho.$
\end{lemma}

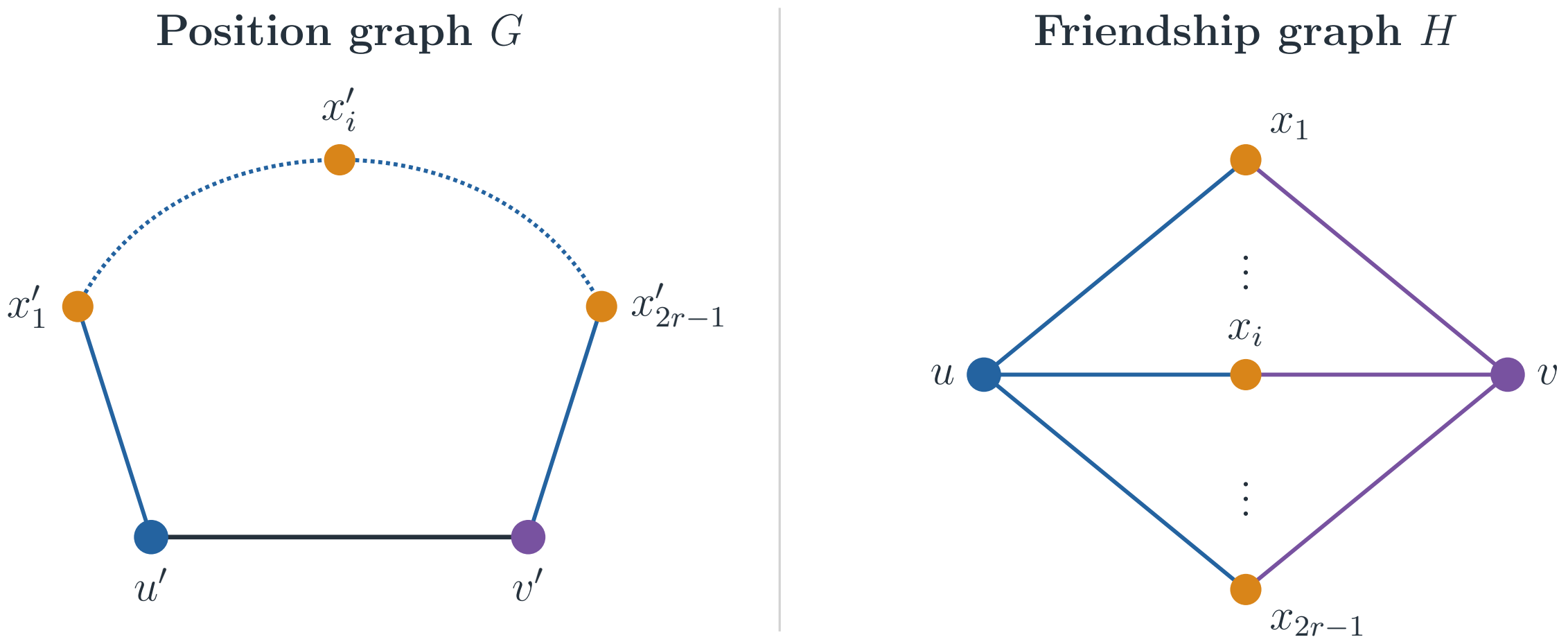
\begin{figure}
\begin{tikzpicture}[
 x=1cm,y=1cm,
 every node/.style={font=\large,text=ink},
 vertex/.style={circle,fill=helperorange,inner sep=3.0pt},
 usite/.style={circle,fill=cycleblue,inner sep=3.3pt},
 vsite/.style={circle,fill=centerpurple,inner sep=3.3pt},
 cycleedge/.style={draw=cycleblue,line width=1.1pt},
 ufriend/.style={draw=cycleblue,line width=1.1pt},
 vfriend/.style={draw=centerpurple,line width=1.1pt},
 lab/.style={inner sep=2pt},
 heading/.style={font=\large\bfseries}
]
% G: the two distinguished positions are adjacent.
% Dotted arcs represent the omitted portions of the cycle.
\node[heading] at (0,3.25) {Position graph $G$};
\coordinate (up) at (-1.80,-1.55);
\coordinate (vp) at (1.80,-1.55);
\coordinate (xfirstp) at (-2.50,0.65);
\coordinate (xmiddlep) at (0,2.05);
\coordinate (xlastp) at (2.50,0.65);
\draw[cycleedge] (up)--(xfirstp);
\draw[cycleedge,densely dotted]
 (xfirstp) .. controls (-2.05,1.55) and (-1.05,2.05) .. (xmiddlep);
\draw[cycleedge,densely dotted]
 (xmiddlep) .. controls (1.05,2.05) and (2.05,1.55) .. (xlastp);
\draw[cycleedge] (xlastp)--(vp);
\draw[draw=ink,line width=1.25pt] (vp)--(up);
\node[usite] at (up) {};
\node[vsite] at (vp) {};
\foreach \pt in {xfirstp,xmiddlep,xlastp}
 \node[vertex] at (\pt) {};
\node[lab,anchor=north] at ($(up)+(0,-0.22)$) {$u'$};
\node[lab,anchor=north] at ($(vp)+(0,-0.22)$) {$v'$};
\node[lab,anchor=east] at ($(xfirstp)+(-0.20,0)$) {$x_1'$};
\node[lab,anchor=south] at ($(xmiddlep)+(0,0.20)$) {$x_i'$};
\node[lab,anchor=west] at ($(xlastp)+(0.20,0)$) {$x_{2r-1}'$};

\draw[gray!35,line width=0.6pt] (4.20,-2.45)--(4.20,3.50);

% H: both distinguished people are friends with every x_i.
% The vertical dots stand for the undisplayed vertices, not edges.
\begin{scope}[xshift=8.65cm]
\node[heading] at (0,3.25) {Friendship graph $H$};
\coordinate (u) at (-2.50,0);
\coordinate (v) at (2.50,0);
\coordinate (xfirst) at (0,2.05);
\coordinate (xmiddle) at (0,0);
\coordinate (xlast) at (0,-2.05);
\foreach \pt in {xfirst,xmiddle,xlast}{
 \draw[ufriend] (u)--(\pt);
 \draw[vfriend] (v)--(\pt);
}
\node[usite] at (u) {};
\node[vsite] at (v) {};
\foreach \pt in {xfirst,xmiddle,xlast}
 \node[vertex] at (\pt) {};
\node[lab,anchor=east] at ($(u)+(-0.20,0)$) {$u$};
\node[lab,anchor=west] at ($(v)+(0.20,0)$) {$v$};
\node[lab,anchor=south west] at ($(xfirst)+(0.15,0.12)$) {$x_1$};
\node[lab,anchor=south] at ($(xmiddle)+(0,0.20)$) {$x_i$};
\node[lab,anchor=north west] at ($(xlast)+(0.15,-0.12)$) {$x_{2r-1}$};
\node at (0,1.08) {$\vdots$};
\node at (0,-1.08) {$\vdots$};
\end{scope}
\end{tikzpicture}
\caption{The graphs $G$ and $H$ in the statement of Lemma~\ref{lem:odd-cycle-exchange}.}
\label{fig:ghconstruction}
\end{figure}

\begin{proof}
Let $G$ have vertices $u',x_1',\ldots,x_{2r-1}',v'$ appearing consecutively in the cycle. Let $H$ have vertices $u, x_1, \ldots, x_{2r-1}, v$ such that $(u,x_i),(v,x_i) \in E(H)$ for every $1\leq i \leq 2r-1.$ Relabeling the remaining vertices according to $\rho$, we may assume $\rho(u')=u$, $\rho(v')=v$, and $\rho(x_i')=x_i$ for every $i$. See Figure~\ref{fig:ghconstruction}.

For each $j=1,\ldots,2r+1$, perform the following swaps:
\[
\begin{cases}
(v\, x_{2r-1}), (v\, x_{2r-2}), \ldots, (v\, x_1)
    & j \text{ odd},\\
(u\, x_{2r-1}), (u\, x_{2r-2}), \ldots, (u\, x_1)
    & j \text{ even}.
\end{cases}
\]
Initially, the first sequence moves $v$ backwards along the
cycle, changing the occupants in the displayed position order from
\[
u,x_1,\ldots,x_{2r-1},v
\quad\text{to}\quad
u,v,x_1,\ldots,x_{2r-1}.
\]
After each step, every $x_i$ advances one position clockwise,
and the consecutive pair $u,v$ advances one position clockwise
with its order reversed. Initially, this pair occupies $v',u'$
in the order $v,u$. After $2r+1$ steps, it again occupies
$v',u'$, now in the order $u,v$, since $2r+1$ is odd.
Meanwhile, each $x_i$ has made one complete circuit and returned
to $x_i'$. Thus the final placement sends $u'$ to $v$,
$v'$ to $u$, and every $x_i'$ to $x_i$, as required.
\end{proof}

\begin{lemma}\label{lem:movability}
Let $r\ge1$, and let $X$ and $Y$ be graphs on $n$ vertices. Suppose that
\begin{itemize}
 \item $X$ is $(2r+1)$-vertex-connected;
 \item every edge of $X$ lies on an odd cycle of length at most $2r+1$;
 \item every two distinct vertices of $Y$ have at least $2r-1$
 common neighbors.
\end{itemize}
Suppose further that, for every placement $\sigma:V(X)\to V(Y)$
and every ordered pair $a,b\in V(Y)$ with
$a'=\sigma^{-1}(a)$, $b'=\sigma^{-1}(b)$, and
$(a',b')\in E(X)$, the movement $a\to b'$
is $(X,Y,2r)$-movable from $\sigma$.
Then $\FS(X,Y)$ is connected.
\end{lemma}
\begin{proof}
Since $X$ is connected, Lemma~\ref{lem:exchangeable} reduces the
claim to exchanging the occupants $u,v$ of any adjacent positions
$u',v'$ from any placement $\sigma$.
Choose an odd cycle through $(u',v')$ with cyclic vertex order
\(
 u',w_1',\ldots,w_{2s-1}',v',
\)
for $1\le s\le r$ and choose distinct common friends $h_1,\ldots,h_{2s-1}$ of $u,v$.

We will now show how to move each $h_j$ onto $w_j'$ for each $j \in \{1,\ldots,2s-1\}.$ We will consider stages $1,\ldots, 2s-1.$ At stage $j$, protect the people
\[
 F_j=\{u,v,h_1,\ldots,h_{j-1}\}.
\]
There are $j+1\le2s\le2r$ protected people.
Their positions can be deleted while leaving $X$ connected, by
$(2r+1)$-vertex-connectivity. Let $\pi_0$ be the current bijection from $V(X)$ to $V(Y)$. We can thus choose a path
\[
v_0'=\pi_0^{-1}(h_j),v_1',\ldots,v_k'=w_j',
\]
avoiding the positions $\pi_0^{-1}(F_j)$.

Initialize $i = 0$ so that $h_j=\pi_i(v_i').$
By $(X,Y,2r)$-movability, there is a sequence of friendly swaps
moving $h_j$ from $v_i'$ to $v_{i+1}'$ without swapping any
person in $F_j$. Denote the resulting placement by $\pi_{i+1}$. We can repeat for $i=0,\ldots,k-1$, and at the end, $h_j$ occupies
$v_k'=w_j'$, and every person in $F_j$ remains fixed.

After completing the above process for every stage, let $\tau$ be the resulting placement. Note that $\tau(u') = u$, $\tau(v') = v$, and $\tau(w_j') = h_j$ for $1 \leq j \leq 2s-1$. Its cycle occupants are
\[
 u,h_1,\ldots,h_{2s-1},v.
\]
Lemma~\ref{lem:odd-cycle-exchange} supplies
friendly swaps on this cycle transforming $\tau$ into
$(u\,v)\circ\tau$ and restoring every $h_j$ to $w_j'$.

We now reverse the entire sequence used to place the common friends.
No swap in that sequence involved $u$ or $v$, so every reverse swap
is still friendly after their exchange. More explicitly, let $P$
be the permutation of people produced by the placement sequence,
so $\tau=P\circ\sigma$. Since $P$ fixes $u,v$, it commutes with
$(u\,v)$, and the final placement is
\[
 P^{-1}\circ(u\,v)\circ P\circ\sigma
 =(u\,v)\circ\sigma.
\]
Thus $u$ and $v$ exchange positions and every other person returns
to their original position. As this holds for every placement and
position edge, $\FS(X,Y)$ is connected. See Figure~\ref{fig:movability-full} for an illustration of these steps when $r = 1.$

\begin{figure}[ht]
\centering
\begin{tikzpicture}[
    scale=0.92,
    transform shape,
    site/.style={
        circle,
        draw,
        thick,
        minimum size=9mm,
        inner sep=0pt,
        font=\small\bfseries,
        fill=white
    },
    empty/.style={
        circle,
        draw,
        thick,
        minimum size=9mm,
        inner sep=0pt,
        fill=white
    },
    vertexlabel/.style={
        font=\scriptsize
    },
    panelcaption/.style={
        font=\small
    },
    translabel/.style={
        font=\small,
        align=center,
        fill=white,
        inner sep=2pt
    },
    arrow/.style={
        ->,
        >=Stealth,
        very thick
    },
    edge/.style={
        thick
    },
    omitted/.style={
        densely dotted,
        thick
    }
]

% ============================================================
% A macro for the underlying graph
%
% #1 = x shift
% #2 = y shift
% #3 = occupant of y'  (empty or y)
% #4 = occupant of x'  (empty or y)
% #5 = occupant of u'
% #6 = occupant of v'
% #7 = panel caption
% ============================================================

% ----------------------------------------------------------------
% (a) INITIAL CONFIGURATION
% ----------------------------------------------------------------
\begin{scope}[xshift=0cm,yshift=5.4cm]

    \coordinate (Y)  at (-2.4,0);
    \coordinate (P1) at (-1.0,0);
    \coordinate (Pk) at (0.8,0);
    \coordinate (X)  at (2.25,0);

    \coordinate (U) at (1.45,1.55);
    \coordinate (V) at (3.05,1.55);

    % path y' --> x'
    \draw[edge]    (Y)--(P1);
    \draw[omitted] (P1)--(Pk);
    \draw[edge]    (Pk)--(X);

    % triangle
    \draw[edge] (U)--(V)--(X)--cycle;

    % occupants
    \node[site]  (Ya)  at (Y)  {$h$};
    \node[empty] (P1a) at (P1) {};
    \node[empty] (Pka) at (Pk) {};
    \node[empty] (Xa)  at (X)  {};

    \node[site] (Ua) at (U) {$u$};
    \node[site] (Va) at (V) {$v$};

    % X-vertex labels
    \node[vertexlabel, below=-1pt of Ya]
        {$h'=v_0'$};
    \node[vertexlabel, below=-1pt of P1a]
        {$v_1'$};
    \node[vertexlabel, below=-1pt of Pka]
        {$v_{k-1}'$};
    \node[vertexlabel, below=-1pt of Xa]
        {$w'=v_k'$};

    \node[vertexlabel, above=-1pt of Ua] {$u'$};
    \node[vertexlabel, above=-1pt of Va] {$v'$};

    % panel caption
    \node[panelcaption] at (0.3,-1.35)
        {\textbf{(a)} Initial configuration $\sigma$};

\end{scope}

% ----------------------------------------------------------------
% (b) y HAS BEEN MOVED TO x'
% ----------------------------------------------------------------
\begin{scope}[xshift=7.7cm,yshift=5.4cm]

    \coordinate (Y)  at (-2.4,0);
    \coordinate (P1) at (-1.0,0);
    \coordinate (Pk) at (0.8,0);
    \coordinate (X)  at (2.25,0);

    \coordinate (U) at (1.45,1.55);
    \coordinate (V) at (3.05,1.55);

    \draw[edge]    (Y)--(P1);
    \draw[omitted] (P1)--(Pk);
    \draw[edge]    (Pk)--(X);

    \draw[edge] (U)--(V)--(X)--cycle;

    \node[empty] (Yb)  at (Y)  {};
    \node[empty] (P1b) at (P1) {};
    \node[empty] (Pkb) at (Pk) {};
    \node[site]  (Xb)  at (X)  {$h$};

    \node[site] (Ub) at (U) {$u$};
    \node[site] (Vb) at (V) {$v$};

    \node[vertexlabel, below=-1pt of Yb]
        {$h'=v_0'$};
    \node[vertexlabel, below=-1pt of P1b]
        {$v_1'$};
    \node[vertexlabel, below=-1pt of Pkb]
        {$v_{k-1}'$};
    \node[vertexlabel, below=-1pt of Xb]
        {$w'=v_k'$};

    \node[vertexlabel, above=-1pt of Ub] {$u'$};
    \node[vertexlabel, above=-1pt of Vb] {$v'$};

    \node[panelcaption] at (0.3,-1.35)
        {\textbf{(b)} After moving $h$ to $w'$: $\tau$};

\end{scope}

% ----------------------------------------------------------------
% TOP TRANSITION
% ----------------------------------------------------------------

\draw[arrow]
    (3.45,6.15) -- (4.85,6.15);

% ----------------------------------------------------------------
% (c) AFTER THE LOCAL THREE-SWAP MANEUVER
% ----------------------------------------------------------------
\begin{scope}[xshift=7.7cm,yshift=0cm]

    \coordinate (Y)  at (-2.4,0);
    \coordinate (P1) at (-1.0,0);
    \coordinate (Pk) at (0.8,0);
    \coordinate (X)  at (2.25,0);

    \coordinate (U) at (1.45,1.55);
    \coordinate (V) at (3.05,1.55);

    \draw[edge]    (Y)--(P1);
    \draw[omitted] (P1)--(Pk);
    \draw[edge]    (Pk)--(X);

    \draw[edge] (U)--(V)--(X)--cycle;

    \node[empty] (Yc)  at (Y)  {};
    \node[empty] (P1c) at (P1) {};
    \node[empty] (Pkc) at (Pk) {};
    \node[site]  (Xc)  at (X)  {$h$};

    % u and v are now exchanged
    \node[site] (Uc) at (U) {$v$};
    \node[site] (Vc) at (V) {$u$};

    \node[vertexlabel, below=-1pt of Yc]
        {$h'=v_0'$};
    \node[vertexlabel, below=-1pt of P1c]
        {$v_1'$};
    \node[vertexlabel, below=-1pt of Pkc]
        {$v_{k-1}'$};
    \node[vertexlabel, below=-1pt of Xc]
        {$w'=v_k'$};

    \node[vertexlabel, above=-1pt of Uc] {$u'$};
    \node[vertexlabel, above=-1pt of Vc] {$v'$};

    \node[panelcaption] at (0.3,-1.35)
        {\textbf{(c)} After the local swaps: $(u\,v) \circ\tau$};

\end{scope}

% ----------------------------------------------------------------
% RIGHT-HAND VERTICAL TRANSITION
% ----------------------------------------------------------------

\draw[arrow]
    (7.7,3.65) -- (7.7,2.55);

% ----------------------------------------------------------------
% (d) FINAL CONFIGURATION
% ----------------------------------------------------------------
\begin{scope}[xshift=0cm,yshift=0cm]

    \coordinate (Y)  at (-2.4,0);
    \coordinate (P1) at (-1.0,0);
    \coordinate (Pk) at (0.8,0);
    \coordinate (X)  at (2.25,0);

    \coordinate (U) at (1.45,1.55);
    \coordinate (V) at (3.05,1.55);

    \draw[edge]    (Y)--(P1);
    \draw[omitted] (P1)--(Pk);
    \draw[edge]    (Pk)--(X);

    \draw[edge] (U)--(V)--(X)--cycle;

    \node[site]  (Yd)  at (Y)  {$h$};
    \node[empty] (P1d) at (P1) {};
    \node[empty] (Pkd) at (Pk) {};
    \node[empty] (Xd)  at (X)  {};

    % u and v remain exchanged
    \node[site] (Ud) at (U) {$v$};
    \node[site] (Vd) at (V) {$u$};

    \node[vertexlabel, below=-1pt of Yd]
        {$h'=v_0'$};
    \node[vertexlabel, below=-1pt of P1d]
        {$v_1'$};
    \node[vertexlabel, below=-1pt of Pkd]
        {$v_{k-1}'$};
    \node[vertexlabel, below=-1pt of Xd]
        {$w'=v_k'$};

    \node[vertexlabel, above=-1pt of Ud] {$u'$};
    \node[vertexlabel, above=-1pt of Vd] {$v'$};

    \node[panelcaption] at (0.3,-1.35)
        {\textbf{(d)} After moving $h$ back to $h'$: $(u\,v) \circ \sigma$};

\end{scope}

% ----------------------------------------------------------------
% BOTTOM TRANSITION
% ----------------------------------------------------------------

\draw[arrow]
    (4.85,0.75) -- (3.45,0.75);

\end{tikzpicture}

\caption{
    The exchangeability argument turning $\sigma$ into $(u\,v)\circ\sigma$. Unlabeled circles contain people whose identities
    are irrelevant to the argument.
}
\label{fig:movability-full}
\end{figure}
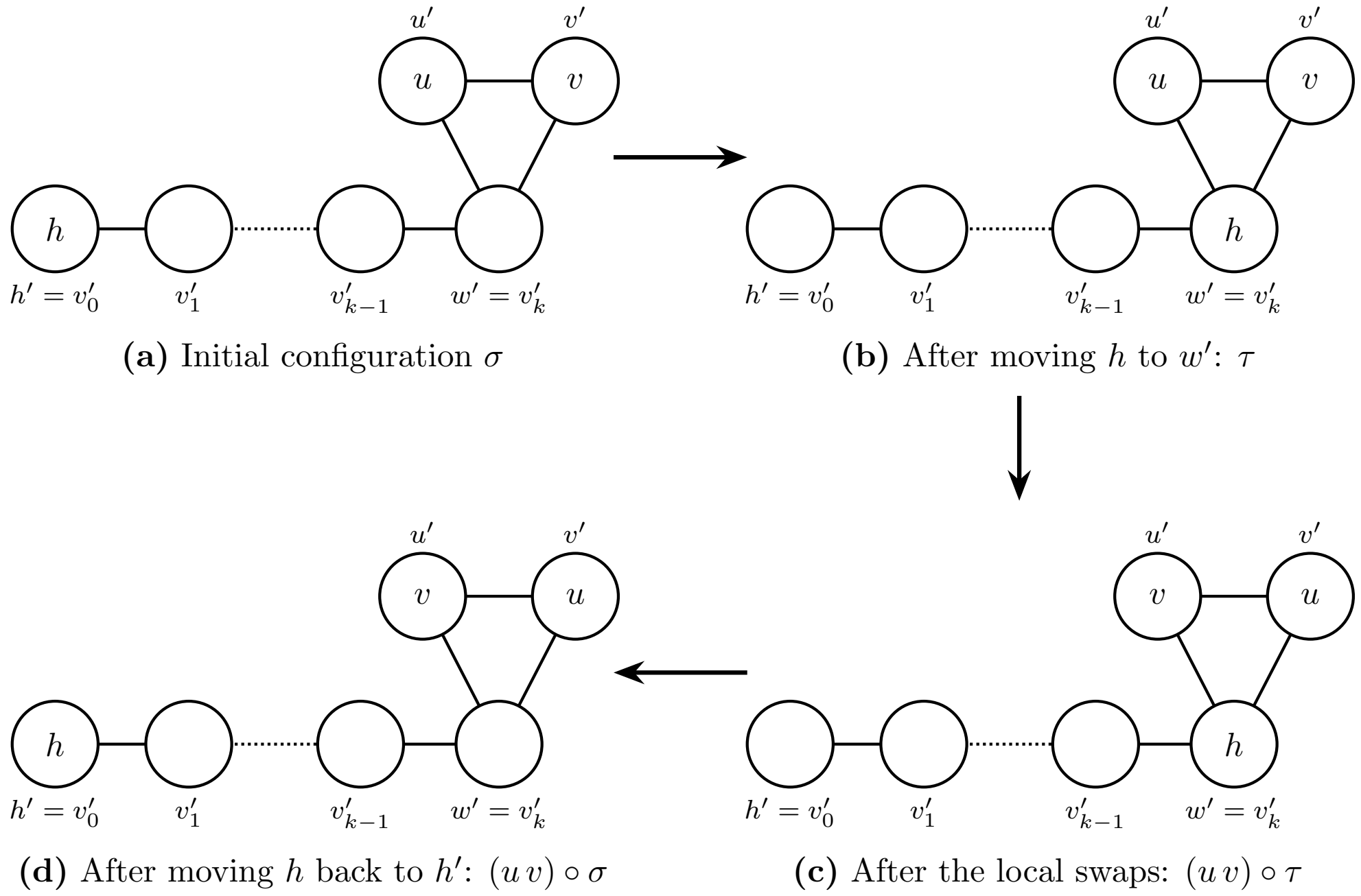
\end{proof}

This type of movement pattern with $r = 1$ (where $r$ defined as in Lemma~\ref{lem:odd-cycle-exchange}) appeared in Bangachev \cite{bangachev2022asymmetric}, e.g., see \cite[Figures 10-12]{bangachev2022asymmetric}, but it was not formalized nor applied to the probabilistic setting where it drastically simplifies constructions.

We can now apply this to the canonical setting using standard random-graph estimates from \cref{sec:preliminaries}. 
\begin{corollary}\label{cor:movability}
Let $X\sim\G(n,p)$ and $Y\sim\G(n,q)$ be independently sampled,
where $p=p(n)\le q=q(n)$ and $pq\ge Cn^{-1}\log n$. Let
\[
r=\left\lceil\frac{2\log n}{\log\log n}\right\rceil.
\]
Then, if with high probability, for every
$\sigma:V(X)\to V(Y)$ and every $a,b\in V(Y)$ where
$a'=\sigma^{-1}(a)$, $b'=\sigma^{-1}(b)$, and $(a',b')\in E(X)$, we have
$a\to b'$ is $(X,Y,2r)$-movable from $\sigma$,
then $\FS(X,Y)$ is connected with high probability.
\end{corollary}

\begin{proof}
We can just apply Lemma~\ref{lem:movability}, noting that $p,q \geq Cn^{-1} \log n.$
Since $2r=o(\log n)$, for sufficiently large $C$,
Lemma~\ref{lem:random-vertex-connectivity} gives that $X$
is $(2r+1)$-vertex-connected,
Lemma~\ref{lem:random-even-paths} gives that every edge
of $X$ is contained in an odd cycle of length at most $2r+1$,
and Lemma~\ref{lem:random-common-neighbors} gives that every
two distinct vertices of $Y$ have at least $2r-1$ common
neighbors, all with high probability.
\end{proof}
The rest of the paper will be devoted to showing $a \to b'$ is $(X,Y,2r)$-movable from $\sigma$ for every $a,b,\sigma$ with high probability.

\section{Construction}\label{sec:construction}
To show movability, we adopt the broad proof strategy of Alon, Defant, and Kravitz \cite{alon2021extremal} which is to find auxiliary configurations in $X$ and $Y$ which allow us to move $a \to b'.$ This is similar to the idea behind Lemma~\ref{lem:movability} where we embedded a construction in order to swap $u$ and $v$.

Unlike Alon, Defant, and Kravitz, we will not need to give a specific graph and embed it into $X$ and $Y.$ Our construction, because of movability, is much simpler, and thus can be found in $X$ and $Y$ by specifying only its structural properties.

% Define an \emph{arm} to be a pair of graphs $G$ and $H$ and a bijection $\rho$ between them. We will let $G$ consist of a cycle $Q$ and a path $P$ extending from a vertex of $Q.$ There are three main distinguished vertices in $G$: $s', \hat s',$ and $c'.$ Vertices $s'$ and $\hat s'$ are distinct vertices on $Q$ where $\hat s'$ is not in $P.$ Vertex $c'$ lies on the endpoint of $P$ not in $Q.$ 
% $H$ on the other hand will be a star centered at $c = \rho(c').$

Fix a placement $\sigma:V(X)\to V(Y)$, distinct people
$a,b\in V(Y)$, and their corresponding positions
$a'=\sigma^{-1}(a)$ and $b'=\sigma^{-1}(b)$, where
$(a',b')\in E(X)$. Also fix a set
$S\subseteq V(Y)\setminus\{a,b\}$ of people we are protecting.

\begin{figure}
\begin{tikzpicture}[
 scale=0.8,x=1cm,y=1cm,
 every node/.style={font=\large,text=ink},
 vertex/.style={circle,fill=ink,inner sep=2.5pt},
 helper/.style={vertex,fill=helperorange,inner sep=3.2pt},
 center/.style={vertex,fill=centerpurple,inner sep=3.2pt},
 port/.style={vertex,fill=cycleblue,inner sep=3.0pt},
 root/.style={vertex,inner sep=3.1pt},
 cycleedge/.style={draw=cycleblue,line width=1.1pt},
 pathedge/.style={draw=centerpurple,line width=1.1pt},
 rootedge/.style={draw=ink,line width=1.1pt},
 friendedge/.style={draw=friendred,line width=1.45pt},
 spoke/.style={draw=centerpurple,line width=0.95pt},
 lab/.style={inner sep=2pt},
 heading/.style={font=\large\bfseries},
 note/.style={font=\normalsize,align=center,text=ink}
]
% Left panel: the required position edges in X.
\node[heading] at (-0.45,4.15) {Position graph $X$};

% First cycle and its pendant path.
\coordinate (s1) at (0,2.05);
\coordinate (t1) at (-0.55,3.05);
\coordinate (u1) at (-1.75,3.05);
\coordinate (l1) at (-2.45,2.05);
\coordinate (e1) at (-1.75,1.05);
\coordinate (d1) at (-0.55,1.05);
\coordinate (c1) at (-4.35,1.05);
\coordinate (p1) at (-3.15,1.05);
\draw[cycleedge] (s1)--(t1)--(u1)--(l1)--(e1)--(d1)--cycle;
\draw[pathedge] (c1)--(p1)--(e1);
\foreach \v in {t1,l1,e1,d1,p1} \node[vertex] at (\v) {};
\node[port] at (s1) {};
\node[helper] at (u1) {};
\node[center] at (c1) {};
\node[lab,anchor=south west] at ($(s1)+(0.12,0.12)$) {$s_1'$};
\node[lab,anchor=south] at ($(u1)+(0,0.20)$) {$\tilde s_1'$};
\node[lab,anchor=north] at ($(c1)+(0,-0.19)$) {$c_1'$};
\node[lab,anchor=south] at (-3.10,1.25) {$P_1$};
\node at (-1.22,2.05) {$Q_1$};

% Second cycle and its pendant path.
\coordinate (s2) at (0,-2.05);
\coordinate (t2) at (-0.55,-3.05);
\coordinate (u2) at (-1.75,-3.05);
\coordinate (l2) at (-2.45,-2.05);
\coordinate (e2) at (-1.75,-1.05);
\coordinate (d2) at (-0.55,-1.05);
\coordinate (c2) at (-4.35,-1.05);
\coordinate (p2) at (-3.15,-1.05);
\draw[cycleedge] (s2)--(t2)--(u2)--(l2)--(e2)--(d2)--cycle;
\draw[pathedge] (c2)--(p2)--(e2);
\foreach \v in {t2,l2,e2,d2,p2} \node[vertex] at (\v) {};
\node[port] at (s2) {};
\node[helper] at (u2) {};
\node[center] at (c2) {};
\node[lab,anchor=north west] at ($(s2)+(0.12,-0.12)$) {$s_2'$};
\node[lab,anchor=north] at ($(u2)+(0,-0.20)$) {$\tilde s_2'$};
\node[lab,anchor=south] at ($(c2)+(0,0.19)$) {$c_2'$};
\node[lab,anchor=north] at (-3.10,-1.25) {$P_2$};
\node at (-1.22,-2.05) {$Q_2$};

% Both ports meet b'; a'b' is the original position edge.
\coordinate (bpos) at (1.85,0);
\coordinate (apos) at (3.45,0);
\draw[rootedge] (s1)--(bpos)--(s2);
\draw[rootedge] (bpos)--(apos);
\node[root] at (bpos) {};
\node[root] at (apos) {};
\node[lab,anchor=north] at ($(bpos)+(0,-0.20)$) {$b'$};
\node[lab,anchor=north] at ($(apos)+(0,-0.20)$) {$a'$};

% Divider between position and friendship graphs.
\draw[gray!35,line width=0.6pt] (4.55,-4.10)--(4.55,4.35);

% Right panel: only the required friendship edges in Y.
% Each star has seven leaves, matching the six cycle vertices and
% one internal pendant-path vertex in its corresponding left arm.
\begin{scope}[xshift=10.1cm]
\node[heading] at (-1.05,4.15) {Friendship graph $Y$};
\coordinate (hc1) at (-2.80,2.35);
\coordinate (hu1) at (0,2.35);
\coordinate (ha) at (2.10,2.35);
\coordinate (hc2) at (-2.80,-2.35);
\coordinate (hu2) at (0,-2.35);
\coordinate (hb) at (2.10,-2.35);

% Six other occupants of the first arm.
\foreach \name/\xx/\yy in {
 A/-4.35/3.50,B/-2.80/3.65,C/-1.25/3.50,
 D/-4.45/2.35,E/-4.35/1.20,F/-1.25/1.20}{
 \coordinate (leaf1\name) at (\xx,\yy);
 \draw[spoke] (hc1)--(leaf1\name);
 \node[vertex] at (leaf1\name) {};
}
\draw[spoke] (hc1)--(hu1);
\node[center] at (hc1) {};
\node[lab,anchor=north] at ($(hc1)+(0,-0.18)$) {$c_1$};

% Six other occupants of the second arm.
\foreach \name/\xx/\yy in {
 A/-4.35/-3.50,B/-2.80/-3.65,C/-1.25/-3.50,
 D/-4.45/-2.35,E/-4.35/-1.20,F/-1.25/-1.20}{
 \coordinate (leaf2\name) at (\xx,\yy);
 \draw[spoke] (hc2)--(leaf2\name);
 \node[vertex] at (leaf2\name) {};
}
\draw[spoke] (hc2)--(hu2);
\node[center] at (hc2) {};
\node[lab,anchor=south] at ($(hc2)+(0,0.18)$) {$c_2$};

% The only required friendships between the arms and a,b.
\draw[friendedge] (ha)--(hu1)--(hu2)--(hb);
\node[helper] at (hu1) {};
\node[helper] at (hu2) {};
\node[root] at (ha) {};
\node[root] at (hb) {};
\node[lab,anchor=south] at ($(hu1)+(0,0.20)$) {$\tilde s_1$};
\node[lab,anchor=north] at ($(hu2)+(0,-0.20)$) {$\tilde s_2$};
\node[lab,anchor=south] at ($(ha)+(0,0.20)$) {$a$};
\node[lab,anchor=north] at ($(hb)+(0,-0.20)$) {$b$};
\end{scope}

\end{tikzpicture}
\caption{The arm construction used in Lemma~\ref{lem:embed->movability}.}
\label{fig:arm}
\end{figure}
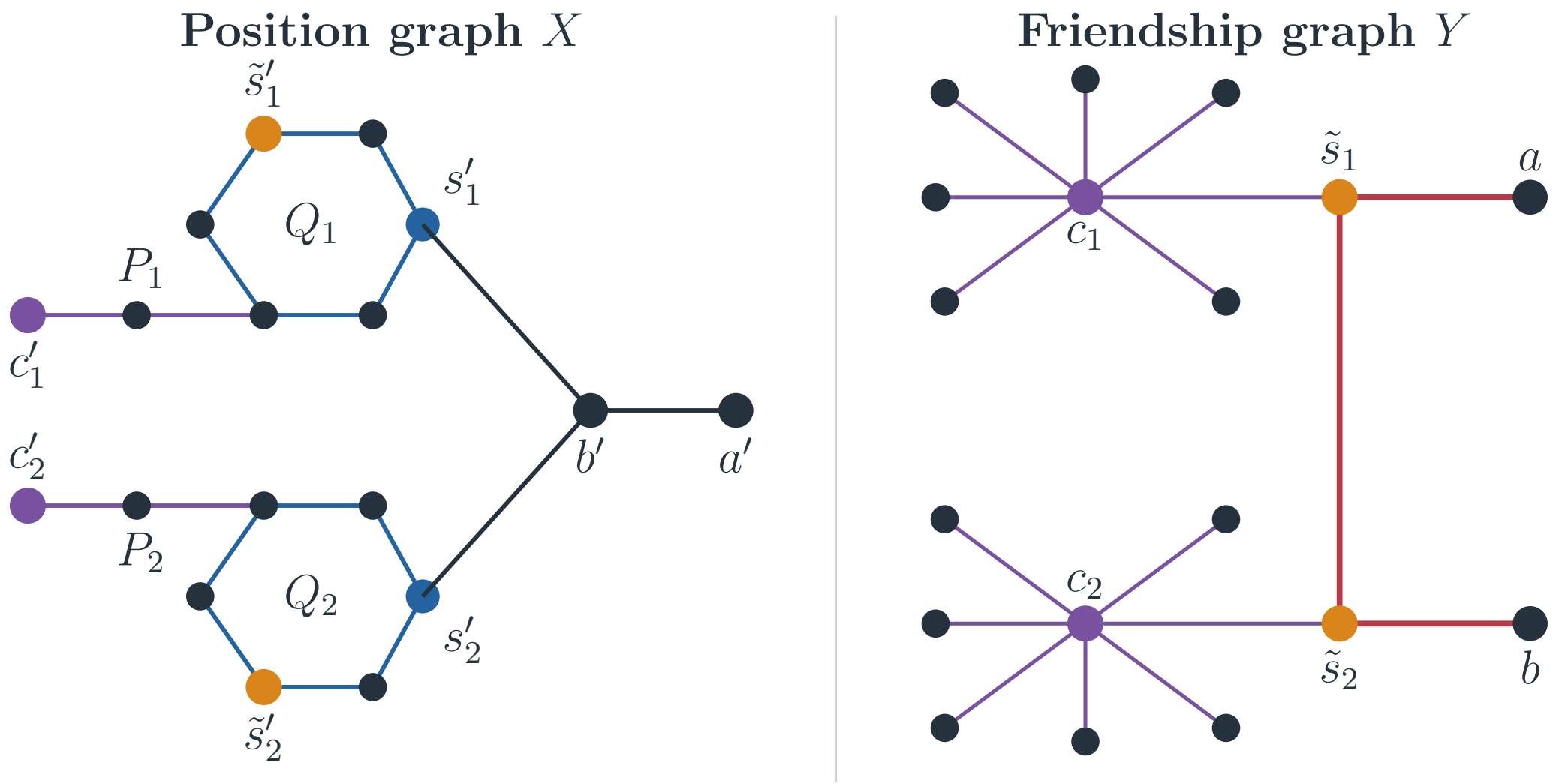

For each $i\in\{1,2\}$, we have the following construction consisting of a simple
cycle $Q_i$ in $X$, together with
\begin{itemize}
    \item a ``port'' position $s_i'\in V(Q_i)$ adjacent to $b'$;
    \item a ``helper'' person $\tilde s_i\in V(Y)$ initially
    at $\tilde s_i'=\sigma^{-1}(\tilde s_i)\in V(Q_i)$;
    \item a ``center'' person $c_i\in V(Y)$ initially at
    $c_i'=\sigma^{-1}(c_i)\notin V(Q_i)$;
\end{itemize}
and a path $P_i$ connecting $c_i'$ to a vertex of
$Q_i$ other than $\tilde s_i'$ such that all internal
vertices of $P_i$ lie outside $Q_i.$ The center $c_i$ is friends with every other person
occupying a position in $V(P_i)\cup V(Q_i)$. We will call this construction an \emph{arm}.
To connect the two arms to $a$ and $b,$ we require
\[
(a,\tilde s_1),\quad
(\tilde s_1,\tilde s_2),\quad
(\tilde s_2,b)\in E(Y).
\]
See Figure~\ref{fig:arm} for an illustration of this construction.

If we can find two such arms which are vertex-disjoint and avoid
$\{a',b'\}\cup\sigma^{-1}(S),$ we can move $a \to b'$ just using friendly swaps within the arms, $a$, and $b.$ Let $G$ be the subgraph of $X$ formed by the two arms, $a',$ and $b',$ and let $H$ be the subgraph of $Y$ formed by the two arms, $a$, and $b.$ Define $\rho = \sigma|_{V(G)}: V(G) \to V(H).$

\begin{lemma}\label{lem:embed->movability}
    There are friendly swaps in $\FS(G,H)$ turning $\rho$ into $\pi$ where $\pi(b') = a.$
\end{lemma}
\begin{proof}
For each $i\in\{1,2\}$, first move $c_i$ along $P_i$ onto
$Q_i$. Every swap is friendly, and $\tilde s_i$ remains on $Q_i$.
We then move $\tilde s_i$ around $Q_i$ to $s_i'$ as follows.
If $\tilde s_i$ is not already at $s_i'$, choose the next
position $v'$ along the cycle towards $s_i'$. Move $c_i$ to $v'$ along the arc of $Q_i$ between $c_i$ and $v'$ which does not contain the current position of $\tilde s_i.$ Finally, swap
$c_i$ with $\tilde s_i$. Repeat until $\tilde s_i$ occupies
$s_i'$. All these swaps are friendly because $c_i$ is friends
with every other person in arm $i$.

Since the arms are disjoint and avoid $a',b'$, these procedures
leave $a,b$ fixed and place $\tilde s_1,\tilde s_2$ at
$s_1',s_2'$, respectively. Now perform the swaps
\[
b\,\tilde s_2,\quad
\tilde s_2\,\tilde s_1,\quad
\tilde s_1\,a.
\]
They are friendly by the three required friendships, and use
the position edges $b's_2'$, $b's_1'$, and $b'a'$, respectively.
The successive occupants of $b'$ are therefore
\[
b,\quad \tilde s_2,\quad \tilde s_1,\quad a.
\]
Thus the resulting placement $\pi$ satisfies $\pi(b')=a$.
\end{proof}

\section{Endgame}
\label{sec:endgame}

In this section, let $X\sim\G(n,p)$ and $Y\sim\G(n,q)$ be
independently sampled, where $p\le q$ and
$pq\ge Cn^{-1}\log n$.
Fix constants $k_1>10$ and $k_2>2k_1$, and choose $C$
sufficiently large afterwards. Throughout, let
\[
r=\left\lceil\frac{2\log n}{\log\log n}\right\rceil,
\qquad
t=\left\lceil k_2q^{-1}\log n\right\rceil.
\]
Note $2r=o(\log n)$ and
$q\ge\sqrt{Cn^{-1}\log n}$.

We will prove that every required movement $a\to b'$ is
$(X,Y,2r)$-movable by embedding the construction using
Lemma~\ref{lem:embed->movability}, and then apply
Corollary~\ref{cor:movability}.

The main idea is to build one arm, using high vertex-connectivity and vertex degrees to construct the cycles and paths necessary. We then build two such
arms in disjoint regions, varying the centers to obtain many possible
helpers $\tilde s_i$. Because of the many choices of arms, we can then choose arms with mutually friendly helpers $\tilde s_1$ and
$\tilde s_2$ among the friends of $a$ and $b$, respectively. This will then construct the two arms necessary.

For fixed $\sigma,a,b,S$ where $|S| \leq 2r$, we call a position \emph{eligible}
if it does not belong to
$\sigma^{-1}(S\cup\{a,b\})$.

\subsection{Finding one arm}
We start by giving simple degree and vertex connectivity conditions for an arm to exist.
\begin{lemma}\label{lem:arm-connectivity}
Choose $c_i\in V(Y)\setminus(S\cup\{a,b\})$, and write
$c_i'=\sigma^{-1}(c_i)$. Let $U$ be any subset of
\[
\sigma^{-1}\bigl(N_Y(c_i)\setminus(S\cup\{a,b\})\bigr).
\]
Suppose that $X[U]$ is $2$-vertex-connected and that both
$b'$ and $c_i'$ have at least two neighbors in $U$.
Then for every $\tilde s_i'\in U$, there is an arm with
$\tilde s_i=\sigma(\tilde s_i')$,
$Q_i\cup P_i\subseteq X[U\cup\{c_i'\}]$, and
$(s_i',b')\in E(X)$.
\end{lemma}

\begin{proof}
Choose $s_i'\in U\cap N_X(b')$ different from
$\tilde s_i'=\sigma^{-1}(\tilde s_i)$.
By $2$-vertex-connectivity and Menger's theorem, two
internally vertex-disjoint paths between $s_i'$ and
$\tilde s_i'$ form a simple cycle $Q_i$.

Choose $z'\in U\cap N_X(c_i')$ different from
$\tilde s_i'$. Since $X[U]-\tilde s_i'$ is connected,
there is a path, possibly of length zero, from $z'$ to
$V(Q_i)\setminus\{\tilde s_i'\}$ avoiding $\tilde s_i'$.
Stop at its first contact with $Q_i$ and prepend the
edge $c_i'z'$. This gives $P_i$, with no internal
vertex in $Q_i$.

Every other person used belongs to $\sigma(U)$ and is
therefore friends with $c_i$, as required.
\end{proof}

\subsection{Admissible partitions}
We now determine the probabilities of the vertex connectivity and degree conditions of Lemma~\ref{lem:arm-connectivity} occurring for a random $U.$

\begin{lemma}\label{lem:sampled-positions}

Fix $V\subseteq V(X)$ of size $n/4+O(1)$, and
independently include each position of $V$ with
probability $q$ to obtain $W$.
With probability at least $1-n^{-4k_1}$ over $X,W$,
the graph $X[W]$ is $(2r+4)$-vertex-connected and
every position of $X$ has at least $2r+4$ neighbors
in $W$.
\end{lemma}

\begin{proof}
Chernoff gives
$nq/8\le M:=|W|\le nq/2$
except with probability $e^{-\Omega(nq)}$.
Conditional on any such $W$, the edges of $X[W]$
remain independent, so
\[
X[W]\sim\G(M,p),
\qquad
Mp\ge npq/8\ge(C/8)\log n.
\]
Since $p\le q$ and $npq\ge C\log n$, we have
$M\ge\sqrt n$ for sufficiently large $n$.
Thus $\log M=\Theta(\log n)$ and
$2r+4\le\log M$.

Lemma~\ref{lem:random-vertex-connectivity},
with $A=8k_1+2$, bounds the failure of
$(2r+4)$-vertex-connectivity by
\[
M^{-(8k_1+2)}\le n^{-4k_1-1},
\]
for sufficiently large $C$.

For each position $v'$, its degree into $W$ has
distribution
\[
\operatorname{Bin}
\bigl(M-\mathbf1_{\{v'\in W\}},p\bigr)
\]
and hence mean $\Omega(npq)$.
Since $2r+4=o(npq)$, Chernoff and a union bound
over all positions give the degree assertion
with failure probability $ne^{-\Omega(npq)}$.
Increasing $C$ makes the total failure probability
at most $n^{-4k_1}$.
\end{proof}

Consider ordered balanced partitions
\[
V(X)=V_1\sqcup H_1\sqcup V_2\sqcup H_2,
\qquad |V_i|,|H_i|=n/4+O(1).
\]
The sets $V_i,H_i$ consist of positions.
For each $i\in\{1,2\}$, let $W_i$ be obtained by
independently retaining each position of $V_i$
with probability $q$.

We call the partition \emph{admissible for $X$} if
\[
B_i(X):=
\mathbb P\bigl(
\text{the conclusion of
Lemma~\ref{lem:sampled-positions}
fails for }W=W_i\mid X
\bigr)
\le n^{-2k_1}
\qquad (i=1,2).
\]

By Lemma~\ref{lem:sampled-positions} and Markov's
inequality, each fixed balanced partition fails
to be admissible for $X$ with probability at most
$2n^{-2k_1}$. Therefore, the expected fraction of
inadmissible partitions is at most $2n^{-2k_1}$.
Another application of Markov shows that, with
probability at least $1-4n^{-2k_1}$, at least half
of all balanced partitions are admissible for $X$.

\subsection{Uniformly finding arms}

\begin{lemma}\label{lem:uniform-arms}
All but fewer than $t$ eligible positions in $V_i$ can serve as
$\tilde s_i'$ in an arm with $c_i'\in H_i$ and all other positions
in $V_i$. For each fixed $X$, this holds simultaneously for every
admissible partition, all choices of $\sigma,a,b,S$ as above, and
both $i\in\{1,2\}$, with probability $1-o(1)$ over $Y$.
\end{lemma}

\begin{proof}
Fix $X$, an admissible partition, $\sigma,a,b,S$,
and $i\in\{1,2\}$.
Let $T$ be any set of $t$ eligible positions in
$V_i$. We show that some position in $T$ can
serve as $\tilde s_i'$.

For each $c_i\notin S\cup\{a,b\}$ with $c_i'=\sigma^{-1}(c_i)\in H_i$, let
\[
W_{c_i}=
\{v'\in V_i:(c_i,\sigma(v'))\in E(Y)\}.
\]
For the fixed placement, $W_{c_i}$ is obtained
by independently retaining each position of
$V_i$ with probability $q$.

By admissibility, $W_{c_i}$ fails the connectivity
or degree conclusion of
Lemma~\ref{lem:sampled-positions}
with probability at most $n^{-2k_1}$.
Moreover,
\[
\mathbb P(W_{c_i}\cap T=\varnothing\mid X)
=(1-q)^t\le e^{-qt}\le n^{-k_2}.
\]

If neither event occurs, choose
$\tilde s_i'\in W_{c_i}\cap T$ and let $U$ consist
of the eligible positions in $W_{c_i}$.
Removing at most $2r+2$ positions preserves
$2$-vertex-connectivity of $X[U]$ and leaves
both $b'$ and $c_i'$ with at least two neighbors
in $U$. Thus Lemma~\ref{lem:arm-connectivity}
gives the required arm.

Consequently, the probability that a fixed
$c_i$ fails to produce an arm with endpoint
in $T$ is at most
\[
n^{-2k_1}+n^{-k_2}\le n^{-k_1}
\]
for sufficiently large $n$.

These events are independent over the choices
of $c_i$, since each depends only on friendship
edges between $c_i$ and people occupying $V_i$,
while the centers occupy the disjoint set $H_i$.
There are at least
$|H_i|-2r-2\ge n/5$ choices of $c_i$.
Hence the probability that every choice fails
is at most
\[
n^{-k_1n/5}.
\]

Taking a union bound over at most $4^n$
partitions, $n!$ placements, $n^2$ choices
of $a,b$, $2^n$ choices each of $S,T$, and
two choices of $i$, the total conditional
failure probability is at most
\[
2\cdot n!\,n^2\,16^n n^{-k_1n/5}=o(1).
\]
Thus every set of $t$ eligible positions in
$V_i$ contains a position supporting an arm,
so fewer than $t$ eligible positions fail
to support one.
\end{proof}

\subsection{Completing the proof}

\begin{theorem}\label{thm:uniform-movability}
With high probability, for every placement
$\sigma$ and every $a,b\in V(Y)$ with
$a'=\sigma^{-1}(a)$, $b'=\sigma^{-1}(b)$, and
$(a',b')\in E(X)$, the movement $a\to b'$ is
$(X,Y,2r)$-movable from $\sigma$.
\end{theorem}

\begin{proof}
With high probability, the following hold:
\begin{enumerate}
    \item $\delta(Y)\ge nq/2$;
    \item every two disjoint sets of $t$ people
    have a friendship between them;
    \item at least half of the balanced partitions
    are admissible for $X$;
    \item the conclusion of
    Lemma~\ref{lem:uniform-arms} holds.
\end{enumerate}
Fix such $X,Y$, and fix $\sigma,a,b,S$ with
$S\subseteq V(Y)\setminus\{a,b\}$ and $|S|\le 2r$.

Since $\delta(Y)\ge nq/2$, a uniformly random
balanced partition satisfies
\[
|\sigma(V_1)\cap N_Y(a)|\ge nq/16,
\qquad
|\sigma(V_2)\cap N_Y(b)|\ge nq/16
\]
with probability $1-o(1)$, by the hypergeometric
Chernoff bound.
Since at least half of the partitions are
admissible for $X$, there exists an admissible
partition satisfying both inequalities.

By Lemma~\ref{lem:uniform-arms}, fewer than $t$
eligible positions in each $V_i$ fail to
support an arm. Moreover, for sufficiently
large $C$,
\[
\frac{nq}{16}-(2r+2)\ge 2t.
\]
Therefore, after excluding the ineligible
positions and those that fail to support arms,
we can choose disjoint sets
\[
A_1\subseteq N_Y(a)\cap\sigma(V_1),
\qquad
A_2\subseteq N_Y(b)\cap\sigma(V_2)
\]
of size $t$, consisting entirely of people occupying
eligible positions that support arms.

By the random-set-edges property, there exist
$\tilde s_1\in A_1$ and $\tilde s_2\in A_2$
such that $(\tilde s_1,\tilde s_2)\in E(Y)$.
Choose the corresponding arms. Since they
lie in the disjoint regions $V_1\cup H_1$
and $V_2\cup H_2$, they are disjoint and
avoid the forbidden positions.

Together with
\[
(a,\tilde s_1),\quad
(\tilde s_1,\tilde s_2),\quad
(\tilde s_2,b)\in E(Y),
\]
Lemma~\ref{lem:embed->movability} gives the
required movement.
\end{proof}

\begin{proof}[Proof of Theorem~\ref{thm:main}]
Inverting placements gives an isomorphism
$\FS(X,Y)\cong\FS(Y,X)$, so we may suppose
$p\le q$. Theorem~\ref{thm:uniform-movability}
and Corollary~\ref{cor:movability} now give
connectivity with high probability.
\end{proof}

\section*{Acknowledgements and AI Disclosure}
The research was conducted with the support of MIT UROP+. The author would like to thank Elchanan Mossel and Rupert Li for helpful advice and discussions. The author would also like to thank Sidhart Krishnan for his support and advice.

The author's contributions were movability, the arm construction, and general proof ideas. ChatGPT Sol 6.1 was used to fill in the details. It also simplified Section~\ref{sec:endgame} as originally the author used a Janson's inequality based argument as in \cite{alon2021extremal}.

The author takes full responsibility for the mathematical content within this paper.

\bibliographystyle{amsinit}
\bibliography{ref}

\end{document}